\documentclass[aos,MSNbibl,dvips]{arximspdf}

\doi{10.1214/13-AOS1160} %kopijuoti is PTS
\volume{41}
\issue{6}
\pubyear{2013}
\firstpage{2768}
\lastpage{2785}

\makeatletter

\newcommand{\rright}{\right}
\newcommand{\lleft}{\left}
\newtheorem{theorem}{Theorem}
\newtheorem{lemma}{Lemma}
\newproclaim{example}{Example}

\makeatother

\begin{document}
\begin{frontmatter}

\title{A complementary set theory for quaternary code~designs}
\runtitle{Quaternary code designs}

\begin{aug}
\author[A]{\fnms{Rahul} \snm{Mukerjee}\ead[label=e1]{rmuk0902@gmail.com}\thanksref{t1}}
\and
\author[B]{\fnms{Boxin} \snm{Tang}\corref{}\ead[label=e2]{boxint@sfu.ca}\thanksref{t2}}
\thankstext{t1}{Supported by the J.~C. Bose National Fellowship of the
Government of
India and a grant from the Indian Institute of Management Calcutta.}
\thankstext{t2}{Supported by the Natural Sciences and Engineering
Research Council of
Canada.}
\runauthor{R. Mukerjee and B. Tang}
\affiliation{Indian Institute of Management Calcutta and Simon Fraser
University}
\address[A]{Indian Institute of Management Calcutta\\
Joka, Diamond Harbour Road\\
Kolkata 700 104\\
India\\
\printead{e1}}

\address[B]{Department of Statistics \\
\quad and Actuarial Science\\
Simon Fraser University\\
Burnaby, BC V5A 1S6\\
Canada\\
\printead{e2}}
\end{aug}

% HISTORY:
\received{\smonth{6} \syear{2013}}

% ABSTRACT
\begin{abstract}
Quaternary code (QC) designs form an attractive class of nonregular
factorial fractions. We develop a complementary set theory for
characterizing optimal QC designs that are highly fractionated in the sense
of accommodating a large number of factors. This is in contrast to existing
theoretical results which work only for a relatively small number of
factors. While the use of imaginary numbers to represent the Gray map
associated with QC designs facilitates the derivation, establishing a link
with foldovers of regular fractions helps in presenting our results in a
neat form.
\end{abstract}

% KEYWORDS
% Pirmas kwd is didziosios raides
%
\begin{keyword}[class=AMS]
\kwd{62K15}
\end{keyword}
\begin{keyword}
\kwd{Foldover}
\kwd{Gray map}
\kwd{highly fractionated design}
\kwd{minimum aberration}
\kwd{minimum moment aberration}
\kwd{projectivity}
\kwd{resolution}
\end{keyword}
\end{frontmatter}

\section{Introduction}\label{sec1} Fractional factorial designs have
received much attention due to their theoretical elegance and practical
applicability to such diverse fields as engineering, agriculture and
medicine. While the literature on regular designs arising from defining
equations is now quite rich, in recent years it has been increasingly
recognized that nonregular designs can potentially perform even better.
See \cite{b6,b12,b15} for detailed reviews and further references.

A significant development in nonregular two-level designs over the last few
years has been the use of quaternary codes (QC) for construction of such
designs, henceforth referred to as QC designs. Xu and Wong \cite{b16} pioneered
work in this direction and this was followed up by \cite{b7,b8,b9,b18}.
As noted
by these authors, QC designs can have an edge over their regular
counterparts under commonly used criteria. Moreover, these designs are
relatively straightforward to construct and have simple design
representation.

The present article aims at developing a theory for optimal QC designs
which accommodate a large number of factors and hence are attractive from
the practical viewpoint of experimental economy. Our theory covers, in
particular, highly fractionated designs with relatively large run sizes.
These have been hitherto unexplored in the context of QC designs but, as
discussed in \cite{b14}, can be of much use in modern day applications. For
instance, our findings are applicable to run sizes 128, 256, 512 and 1024
with 96--112, 224--240, 448--480 and 960--992 factors, respectively. No result,
either computational or theoretical, is as yet available on optimal QC
designs in these cases. For example, the 128- or 256-run design tables in
\cite{b16} cover up to 64 factors. On the other hand, the results in
\cite{b7,b9,b18}
on $1/4$th, $1/8$th, $1/16$th or $1/64$th fractions, though theoretically
appealing, are applicable only when the number of factors is small compared
to the run size, for example, they together cater only to the cases of 10,
11, 12 or 14 factors for run size 256.

Indeed, the existing approaches for theoretical study of QC designs, such
as those based on induction \cite{b9}, trigonometric formulation \cite
{b18} or code
arithmetic~\cite{b7}, get increasingly involved and unmanageable as the
degree of
fractionation increases. To overcome this difficulty, we develop a
complementary set theory, via the use of imaginary numbers, which
considerably facilitates the task of finding optimal QC designs in such
situations. While the use of complementary sets for the study of QC designs
is inspired by the corresponding development in the regular case
\cite{b2,b3,b11}, neither our final results nor their method of
derivation can be
anticipated from the latter.

The commonly used optimality criteria in selecting factorial fractions are
resolution, aberration and projectivity, as introduced briefly later in
this section. In our setup where the number of factors exceeds half the run
size, regular designs have~resolution three and projectivity two, while
following \cite{b5,b16}, QC designs have resolution at least 3.5 and
projectivity
at least three. In other words, with regard to resolution and projectivity,
QC designs have an edge over regular ones. So, there will be a strong case
in favor of QC designs, justifying their use in practice, provided they
compete well with regular designs under the minimum aberration (MA)
criterion as well. This necessitates identification, in our highly
fractionated setup, of minimum aberration quaternary code (MA QC) designs,
which is precisely the focus of the present work. Thus, in addition to
strengthening the currently available theory of QC designs, our results
facilitate their comparison with regular MA designs and help in making a
choice between the two.

Before concluding the \hyperref[sec1]{Introduction}, we briefly recall some definitions.
A~two-level design $D$ in $N$ runs and $q$ factors is represented by an $N
\times q$ matrix with elements $\pm1$, where the rows and columns are
identified with the runs and factors, respectively. The \textit{aliasing
index} of any subset $H$ of $k$ columns of $D$ is defined as $\rho_{k}(H;D)
=|\mathrm{mean}\{\mathrm{Schur}(H)\}|$. Here Schur($H$) is the Schur product of the columns
in $H$, {that is}, each element of Schur($H$) is the product of the
corresponding elements of the columns in $H$, and mean $\{\mathrm{Schur}(H)\}$ is the
arithmetic mean of the elements of Schur($H$); see \cite{b4,b5}.
Clearly, we have
$0 \le\rho_{k}(H;D) \le1$. The columns in $S$ are fully aliased,
partially aliased or unaliased according as $\rho_{k}(H;D)$ equals 1, lies
strictly between 0 and 1, or equals 0, respectively. For $1 \le k \le q$,
let
\begin{equation}
\label{eq1} \rho_{k,\max} (D)=\max\rho_{k}(H;D),\qquad
A_{k}(D)= \sum \bigl\{ \rho_{k}(H;D) \bigr
\}^{2},
\end{equation}
the maximum and the sum being over all $H$ of cardinality $k$. Write $r$
for the smallest integer such that $\rho_{r,\max} (D)> 0$. Then the
\textit{resolution} of $D$ is defined \cite{b5} as $R(D)=r + 1 - \rho
_{r,\max}
(D)$, while the MA criterion calls for sequential minimization of the
components of the vector ($A_{1}(D),\ldots, A_{q}(D)$), known as the
\textit{wordlength pattern} (WLP) of $D$ \cite{b10,b17}. These
definitions apply
to all two-level designs, regular or not, and reduce to the corresponding
combinatorial definitions in the regular case. Finally, following \cite{b1},
design $D$ is said to have \textit{projectivity} $p$ if $p$ is the largest
integer such that every $p$-factor projection of $D$ contains a complete
$2^{p}$ factorial design, possibly with some points replicated.

\section{\texorpdfstring{Quaternary code designs in $2^{2n}$ runs and an even number of factors}
{Quaternary code designs in 2 2n runs and an even number of factors}}\label{sec2}

\subsection{Preliminaries}\label{sec2.1} Let $C$ be the quaternary
linear code given by the $n \times s$ generator matrix $G = [g_{1}\cdots
g_{s}]$ whose columns are $n \times1$ vectors, $n \ge2$, over the set
of integers $Z_{4} = \{ 0,1,2,3\}$ (mod 4). The code $C$, consisting of
$4^{n}\ (=2^{2n})$ codewords, each of size $s$, can be described as
\begin{equation}
\label{eq2} C = \bigl\{ \bigl(u'g_{1},
\ldots,u'g_{s} \bigr)\dvtx u = (u_{1},
\ldots,u_{n})',u_{j} \in Z_{4},1 \le
j \le n \bigr\},
\end{equation}
where the primes stand for transposition and each of
$u'g_{1},\ldots,u'g_{s}$ is reduced mod~4. The \textit{Gray map}, which
replaces each element of $Z_{4}$ with a pair of two symbols according to the
rule
\begin{equation}
\label{eq3}0\to(1,1),\qquad 1\to (1, -1),\qquad 2\to (-1, -1),\qquad 3\to
(-1,1),
\end{equation}
transforms $C$ into a binary code $D$, called the binary image of $C$. With
its codewords as rows, $D$ is a $2^{2n} \times(2s)$ matrix having entries
$\pm1$. In this sense, with columns and rows identified with factors and
runs, respectively, $D$ represents a QC design in $2s$ two-level factors
and $2^{2n}$ runs.

In order to meet the essential design objective of keeping all main effects
orthogonally estimable at least when interactions are absent, the QC design
$D$ constructed as above must be an orthogonal array of strength two. Xu
and Wong \cite{b16} showed that $D$ meets this basic requirement if and
only if
(a) none of $g_{1},\ldots,g_{s}$ has each entry even, that is, 0 or 2, and (b)
no two of these vectors are multiples of each other over $Z_{4}$. They also
noted that if these conditions hold, then $D$ has resolution at least 3.5
and hence projectivity at least three. So, hereafter we consider only those
QC designs which satisfy (a) and (b). Then each $g_{1},\ldots,g_{s}$ has an
odd element and, without loss of generality, the first odd element in each
of these is 1. Following \cite{b16} again, there are altogether $v =
(4^{n} -
2^{n})/2$ such vectors over $Z_{4}$. Let $\Omega$ be the collection of
these $v$ vectors. For example, if $n = 2$, then $\Omega=\{
(0,1)',(1,0)',(1,1)',(1,2)',(1,3)',(2,1)'\}$. The task of obtaining an MA
QC design then amounts to finding a subset $S = \{ g_{1},\ldots,g_{s}\}$ of
$\Omega$ such that the design $D$ arising from $G = [g_{1}\cdots g_{s}]$
sequentially minimizes $A_{3}(D),\ldots, A_{q}(D)$, where $q = 2s$. The
first two terms of the WLP, $A_{1}(D)$ and $A_{2}(D)$, are dropped here
since they equal zero for any such $D$. The case $s = v$ is trivial, for
then $S =\Omega$ is the unique choice of $S$. Therefore, in the rest of
this section, we consider $s < v$.

For highly fractionated designs which form our main focus, $s$ is large.
Hence, instead of considering $S = \{ g_{1},\ldots,g_{s}\}$ directly, it will
be more convenient and insightful to work with $\bar{S}$, the
complement of
$S$ in $\Omega$. This is motivated by what one does in the study of regular
designs, but there are major differences. For example, the counterpart of
$\Omega$ in the regular case is a finite projective geometry, while the
vectors in our $\Omega$ are not even on a finite field. This warrants the
development of new techniques of proof for the present problem.

\subsection{MA criterion in terms of complementary set}\label{sec2.2}
We now proceed to formulate the MA criterion in terms of the
complementary set $\bar{S}$ introduced above. To that effect, from
(\ref{eq2}) note that the $2^{2n}$ rows of $D$ can be indexed by the
$2^{2n}$ vectors $u = (u_{1},\ldots,u_{n})'$ over $Z_{4}$. Let $\Delta$
be the collection of all such $u$. For any $u\ ( \in\Delta)$, write
$\theta'_{u}$ for the corresponding row of $D$. Clearly, as $D$ has
$2s$ columns and elements $\pm1$, $\theta '_{u}\theta_{w} = 2c_{uw} -
2s$, where $c_{uw}$ is the number of coincidences between
$\theta'_{u}$ and $\theta'_{w}$. In view of the equivalence between the
MA and minimum moment aberration (MMA) criteria as established in
\cite{b13}, it is clear that sequential minimization of
$A_{3}(D),\ldots, A_{q}(D)$ is equivalent to that of
$M_{3}(D),\ldots,M_{q}(D)$, where
\begin{equation}
\label{eq4} M_{k}(D)= \sum_{u \in\Delta} \sum
_{w \in\Delta} \bigl(\theta '_{u}
\theta_{w} \bigr)^{k}, \qquad 3 \le k \le q.
\end{equation}

In order to achieve further simplification, taking due cognizance of the
structure of a QC design, let $i = \sqrt{ - 1}$ and, for any integer $z$,
write
\begin{equation}
\label{eq5} \psi(z) = \bigl(i^{z} + i^{3 - z} \bigr)/(1 - i).
\end{equation}
Since the Gray map (\ref{eq3}) is equivalent to $z \to(\psi( -
z),\psi(z))$, $z
\in Z_{4}$, by (\ref{eq2}), the row $\theta'_{u}$ of $D$ can be
written explicitly
as
\begin{equation}
\label{eq6} \theta'_{u}= \bigl(\psi \bigl( -
u'g_{1} \bigr),\psi \bigl(u'g_{1}
\bigr),\ldots,\psi \bigl( - u'g_{s} \bigr),\psi
\bigl(u'g_{s} \bigr) \bigr).
\end{equation}
For any $u \in\Delta$, let $\sigma_{u}$ be the sum of the elements of
$\theta'_{u}$. Then, from (\ref{eq5}) and (\ref{eq6}), after some algebra,
\begin{equation}
\label{eq7} \theta'_{u}\theta_{w} =
\sigma_{u - w},\qquad u,w \in\Delta,
\end{equation}
where $u - w$ is reduced mod 4; cf. Theorem 3 in \cite{b16}. Now, for
any fixed
$w \in\Delta$, as $u$ equals every member of $\Delta$, so does $u - w$.
Hence, by (\ref{eq4}) and (\ref{eq7}), $M_{k}(D)=2^{2n}m_{k}(D)$, where
\begin{equation}
\label{eq8} m_{k}(D)=\sum_{u \in\Delta}
\sigma_{u}^{k}, \qquad 3 \le k \le q,
\end{equation}
and sequential minimization of $A_{3}(D),\ldots, A_{q}(D)$ reduces to that of
$m_{3}(D),\break\ldots, m_{q}(D)$. Note that the $m_{k}(D)$ involve only the row
totals of $D$ rather than scalar products of rows, and hence are much
simpler than the $M_{k}(D)$.

We next express the $\sigma_{u}$ and hence the quantities $m_{k}(D)$ in terms
of the complementary set $\bar{S}$. Recall that $S = \{
g_{1},\ldots,g_{s}\}$.
Hence, if we write $\Omega=\{ g_{1},\ldots,g_{v}\}$, where $v = (4^{n} -
2^{n})/2$ as before, then $\bar{S} = \{ g_{s + 1},\ldots,g_{v}\}$. Analogously
to $D$, let $\bar{D}$ be the QC design, in $2(v - s)$ two-level factors and
$2^{2n}$ runs, arising from the generator matrix $\bar{G} = [g_{s +
1}\cdots g_{v}]$. For $u \in\Delta$, write $\bar{\sigma}_{u}$ as the sum of
elements of the row of $\bar{D}$ which is indexed by $u$. By (\ref
{eq5}) and the
counterpart of (\ref{eq6}) for $\bar{D}$,
\begin{equation}
\label{eq9} \bar{\sigma}_{u}=\sum_{j = s + 1}^{v}
\bigl(i^{u'g_{j}} + i^{ -
u'g_{j}} \bigr),\qquad u \in\Delta.
\end{equation}
Also, let $\Delta_{0}$ consist of the $2^{n}$ vectors in $\Delta$ which
have all elements even, and define $\delta_{u}$ as 1 or 0 according to
whether $u$ belongs to $\Delta_{0}$ or not. Trivially, by (\ref{eq5}),
(\ref{eq6}) and
(\ref{eq9}),
\begin{equation}
\label{eq10} \sigma_{u}=2s,\qquad \bar{\sigma}_{u}=2(v -
s)\qquad \mbox{if }u = (0,\ldots,0)'.
\end{equation}

The following lemma connects $\sigma_{u}$ and $\bar{\sigma}_{u}$ for
nonnull $u$.

\begin{lemma}\label{le1}
Let $u( \in\Delta)$ be nonnull. Then
$\sigma_{u}=- (2^{n}\delta_{u} + \bar{\sigma}_{u})$.
\end{lemma}

By (\ref{eq8}), (\ref{eq10}) and Lemma \ref{le1}, $m_{k}(D)= \mathrm{constant} +( -
1)^{k}\bar{m}_{k}$ for
each $k$, where
\begin{equation}
\label{eq11} \bar{m}_{k} = \sum_{u \in\Delta}
\bigl(2^{n}\delta_{u} + \bar{\sigma}_{u}
\bigr)^{k},
\end{equation}
and the constant does not depend on $D$. Hence, in the quest of an MA QC
design, one needs to find $\bar{S}$ so as to maximize $\bar{m}_{3}$, then
minimize $\bar{m}_{4}$, then maximize $\bar{m}_{5}$, and so on.

\subsection{Characterization of MA designs}\label{sec2.3}
Since each $g_{j}$ has an odd element, from~(\ref{eq9}), (\ref{eq11})
and the definitions of $\Delta_{0}$ and $\delta_{u}$, arguments similar
to but simpler than those in the proof of Theorem\vadjust{\goodbreak} \ref{th1} below show
that $\bar{m}_{3}= \mathrm{constant} +F_{3}$, where the constant does not
depend on $\bar{S}$ and
\begin{equation}
\label{eq12} F_{3} = 3 \bigl(2^{n} \bigr) \sum
_{u \in
\Delta_{0}}\bar{\sigma}_{u}^{2}+ \sum
_{u \in\Delta} \bar{\sigma}_{u}^{3}.
\end{equation}

We first explore $\bar{S}$ so as to maximize $\bar{m}_{3}$ or,
equivalently, $F_{3}$. The set $\bar{S} = \{ g_{s + 1},\ldots,g_{v}\}$ is
called \textit{even} if $g_{j} + g_{k}$ has all elements even for every $s
+ 1 \le j,k \le v$. Then the following theorem, which is a main result of
this section, holds.

\begin{theorem}\label{th1}
\textup{(a)} The inequality $F_{3} \le3(2^{2n + 2})(v - s)^{2}$
holds for every~$\bar{S}$.

\textup{(b)} Equality holds in \textup{(a)} if and only if $\bar{S}$
is even.
\end{theorem}

In order to apply Theorem \ref{th1}, one needs to know when an even set $\bar{S}$
exists. The next lemma settles this issue.

\begin{lemma}\label{le2}
There exists an even set $\bar{S}$ if and only
if $v - s \le2^{n - 1}$.
\end{lemma}

For $v - s \le2^{n - 1}$, Theorem \ref{th1} and Lemma \ref{le2}
significantly reduce the
problem of finding an MA QC design, since one needs to consider only sets
$\bar{S} = \{ g_{s + 1},\ldots,g_{v}\}$, which are even. Then each row of
$\bar{G} = [g_{s + 1}\cdots g_{v}]$ has either all elements odd or all elements
even. Indeed, $\bar{G}$ has a row with all elements odd, because $\bar{S}
\subset\Omega$. Since the row space of $\bar{G}$ remains invariant under
elementary row operations, without loss of generality, let
\begin{equation}
\label{eq13} \bar{G}=\lleft[ \matrix{ 1& 1'_{v - s - 1}
\cr
0_{n - 1}&2B } \rright],
\end{equation}
where $0_{n - 1}$ is the column vector of $n - 1$ zeros, $1'_{v - s -
1}$ is
the row vector of $v - s - 1$ ones, and $B$ is an $(n - 1) \times(v - s -
1)$ binary matrix. Clearly, $B$ can be interpreted as the generator matrix
of a regular design, say, $d$, in $v - s - 1$ two-level factors and
$2^{n -1}$ runs (these runs are not distinct in case $B$ has less than full row
rank). Our next theorem characterizes the MA property of $\bar{S}$ in terms
of the much simpler regular design $d$. Here $A_{k}(d)$, $k \ge1$, denotes
the WLP of $d$ and $A{}_{1}(d) = A_{2}(d) = 0$, because the columns of
$\bar{G}$ are distinct and, as a result, those of $B$ are distinct and
nonnull.

\begin{theorem}\label{th2}
Let $v - s \le2^{n - 1}$. Then
$\bar{S}$ yields an MA QC design in 2s factors and
$2^{2n}$ runs if and only if the matrix B is so chosen that the
associated two-level regular design d sequentially minimizes $A_{2r -
1}(d) + A_{2r}(d)$, for $r = 2,3,\ldots,$ etc.
\end{theorem}

Theorem \ref{th2} suggests a connection with the full foldover of $d$, a point
which is confirmed by its proof in the \hyperref[app]{Appendix}. To
apply this theorem, one
needs to consider all nonisomorphic choices of $d$ and select one from
among them meeting the condition of the theorem. The extensive tables of
two-level regular designs available in the literature are very useful in
this regard. We now present two illustrative examples followed by a design
table where, for notational simplicity, any binary column vector with 1 in
positions $h_{1},\ldots,h_{r}$ and zeros elsewhere is represented by
$h_{1}\cdots h_{r}$. Thus, the binary matrix, with columns
$(1,0,0)',(0,1,0)',(1,1,0)'$ and $(0,0,1)'$ is denoted by [1 2 12 3].

\begin{example}\label{ex1}
Let $n = 4$, that is, $v = 120$, and $v - s= 5$. Then there are
two nonisomorphic choices of $B$, namely, [1 2 12 3] and [1 2 3 123]. Both
entail $A_{3}(d) + A_{4}(d) = 1$ and hence satisfy the condition of
Theorem~\ref{th2}, leading to MA QC designs in 230 factors and 256 runs.
\end{example}

\begin{example}\label{ex2} Let $n = 5$, that is,
$v = 496$, and $v - s= 10$. Table 3A.2 in
\cite{b6} lists all nonisomorphic choices of $B$ along with the WLPs
($A_{3}(d),\ldots,A_{9}(d)$) of the corresponding two-level regular designs
$d$. These are as shown below:
\begin{longlist}[(iii)]
\item[(i)] $B = {}$[1 2 12 3 13 4 14 234 1234], WLP${}={}$(4, 14, 8, 0, 4, 1, 0);

\item[(ii)] $B = {}$[1 2 12 3 13 4 24 34 1234], WLP${}={}$(6, 9, 9, 6, 0, 0, 1);

\item[(iii)] $B = {}$[1 2 12 3 13 23 4 14 234], WLP${}={}$(6, 10, 8, 4, 2, 1, 0);

\item[(iv)] $B = {}$[1 2 12 3 13 23 4 14 24], WLP${}={}$(7, 9, 6, 6, 3, 0, 0);

\item[(v)] $B = {}$[1 2 12 3 13 23 123 4 14], WLP${}={}$(8, 10, 4, 4, 4, 1, 0).
\end{longlist}
The choice in (ii) uniquely minimizes $A_{3}(d) + A_{4}(d)$ and hence yields
an MA QC design in 972 factors and 1024 runs. Incidentally, the design $d$
associated with~(ii) does not itself have MA as a two-level regular
fraction in $9\ (= v - s - 1)$ factors and $16\ ( = 2^{n - 1})$ runs.
\end{example}

For each $n = 3,4,5$ and $v - s \le2^{n - 1}$, Table~\ref{tab1} shows the columns
of $B$ so that the associated $d$ satisfies the condition of Theorem \ref{th2} and
hence yields an MA QC design in $2s$ factors and $2^{2n}$ runs. In the
event of nonuniqueness as in Example \ref{ex1}, only one such $B$ is shown. The
cases $v - s= 0$ and 1 are not considered in Table~\ref{tab1} because if $v - s= 0$,
then $\bar{S} = \Omega$ is the unique choice of $\bar{S}$, while if $v
- s=
1$, then the matrix $B$ does not arise in (\ref{eq13}) and it suffices
to take
$\bar{S}$ as the singleton set consisting of $(1,0,\ldots,0)'$.

\begin{table}
\caption{Choice of B leading to an MA QC design via Theorem \protect\ref{th2}}\label{tab1}
\begin{tabular*}{\textwidth}{@{\extracolsep{\fill}}llll@{}}
\hline
$n = 3$& &$n = 5$&\\
$v - s$& Columns of $B$&$v - s$& Columns of $B$\\
2& 1& \phantom{0}2& 1\\
3 &1\ 2 & \phantom{0}3& 1\ 2\\
4 & 1\ 2 12& \phantom{0}4& 1\ 2\ 3\\
&& \phantom{0}5& 1\ 2\ 3\ 4\\
$n = 4$ & &\phantom{0}6& 1\ 2\ 3\ 4\ 1234\\
$v - s$&  Columns of $B$ &\phantom{0}7& 1\ 2\ 12\ 3\ 4\ 34\\
2& 1& \phantom{0}8& 1\ 2\ 12\ 3\ 13\ 4\ 24\\
3 & 1\ 2& \phantom{0}9& 1\ 2\ 12\ 3\ 13\ 4\ 24\ 34\\
4& 1\ 2\ 3& 10& 1\ 2\ 12\ 3\ 13\ 4\ 24\ 34\ 1234\\
5 & 1\ 2\ 12\ 3& 11& 1\ 2\ 12\ 3\ 13\ 23\ 4\ 14\ 24\ 34\\
6& 1\ 2\ 12\ 3\ 13& 12& 1\ 2\ 12\ 3\ 13\ 23\ 123\ 4\ 14\ 24\ 34\\
7& 1\ 2\ 12\ 3\ 13\ 23& 13& 1\ 2\ 12\ 3\ 13\ 23\ 123\ 4\ 14\ 24\ 124\
34\\
8& 1\ 2\ 12\ 3\ 13\ 23\ 123& 14& 1\ 2\ 12\ 3\ 13\ 23\ 123\ 4\ 14\ 24\ 124\ 34\
134\\
&& 15& 1\ 2\ 12\ 3\ 13\ 23\ 123\ 4\ 14\ 24\ 124\ 34\ 134\ 234\\
&& 16& 1\ 2\ 12\ 3\ 13\ 23\ 123\ 4\
14\ 24\ 124\ 34\ 134\ 234\ 1234\\
\hline
\end{tabular*}
\end{table}

\section{\texorpdfstring{Quaternary code designs in $2^{2n}$ runs and an odd number of factors}
{Quaternary code designs in 2 2n runs and an odd number of factors}}\label{sec3}

A~QC design with an odd number of
factors is
constructed as follows. First construct a QC design $D_{\mathrm{even}}$ in
$2(s + 1)$ two-level factors and $2^{2n}$ runs as in Section~\ref{sec2} starting from
the generator matrix [$g_{1}\cdots g_{s}\ g_{s + 1}$], where $s + 1 \le v$. Each
vector $g_{1},\ldots,g_{s},g_{s + 1}$ contributes two columns to
$D_{\mathrm{even}}$ via the Gray map (\ref{eq3}). Delete a column of
$D_{\mathrm{even}}$ to get a QC design $D_{\mathrm{odd}}$ in $2s + 1$ factors
and $2^{2n}$ runs. Without loss of generality, suppose the second column
contributed by $g_{s + 1}$ is deleted. Then as in (\ref{eq6}), the rows of
$D_{\mathrm{odd}}$ are given by
\begin{eqnarray}
\label{eq14} \theta'_{\mathrm{odd},u}= \bigl(\psi \bigl( -
u'g_{1} \bigr),\psi \bigl(u'g_{1}
\bigr),\ldots,\psi \bigl( - u'g_{s} \bigr),\psi
\bigl(u'g_{s} \bigr),\psi \bigl( - u'g_{s + 1}
\bigr) \bigr),
\nonumber
\\[-8pt]
\\[-8pt]
\eqntext{u \in\Delta.}
\end{eqnarray}
Let $\sigma_{\mathrm{odd},u}$ be the sum of the elements of $\theta
'_{\mathrm{odd},u}$. By (\ref{eq5}) and (\ref{eq14}), analogously to~(\ref{eq7}), for $u,w \in
\Delta$, the scalar product $\theta
'_{\mathrm{odd},u}\theta_{\mathrm{odd},w}$ equals $\sigma_{\mathrm
{odd},u -
w}$ or $\sigma_{\mathrm{odd},w - u}$ according to whether $w'g_{s + 1}$ is
even or odd, respectively. Here $u - w$ and $w - u$ are reduced mod 4.
Hence, arguing as in Section~\ref{sec2}, finding an MA QC design calls for
sequential minimization of $m_{k}(D_{\mathrm{odd}})$, $3 \le k \le q$,
where
\begin{equation}
\label{eq15} m_{k}(D_{\mathrm{odd}})=\sum
_{u \in\Delta} \sigma_{\mathrm{odd},u}^{k}.
\end{equation}

It again helps to consider the set $\bar{S} = \{ g_{s + 1},\ldots,g_{v}\}$
even though it is no longer a truly complementary set because of the
partial contribution of $g_{s + 1}$ to $D_{\mathrm{odd}}$. Let
\begin{equation}
\label{eq16} \bar{\sigma}_{\mathrm{odd},u}=\psi \bigl(u'g_{s + 1}
\bigr)+\sum_{j
= s +
2}^{v} \bigl\{ \psi \bigl( -
u'g_{j} \bigr) + \psi \bigl(u'g_{j}
\bigr) \bigr\},\qquad u \in\Delta.
\end{equation}
Evidently, $\sigma_{\mathrm{odd},u}= \sigma_{u} + \psi( - u'g_{s + 1})$
and $\bar{\sigma}_{\mathrm{odd},u}=\bar{\sigma}_{u} - \psi( - u'g_{s +
1})$ for each~$u$, recalling the definitions of $\sigma_{u}$ and
$\bar{\sigma}_{u}$. Hence, by (\ref{eq5}) and (\ref{eq10}), if $u =
(0,\ldots,0)'$, then
$\sigma_{\mathrm{odd},u}=2s + 1$ and $\bar{\sigma}_{\mathrm{odd},u}=
2(v -
s) - 1$, while by Lemma \ref{le1}, if $u( \in\Delta)$ is nonnull, then
$\sigma_{\mathrm{odd},u}= - (2^{n}\delta_{u} +
\bar{\sigma}_{\mathrm{odd},u})$. As a result, by (\ref{eq15}),
$m_{k}(D_{\mathrm{odd}})= \mathrm{constant} +( - 1)^{k}\bar{m}_{\mathrm{odd},k}$
for each $k$, where
\begin{equation}
\label{eq17} \bar{m}_{\mathrm{odd},k} = \sum_{u \in\Delta}
\bigl(2^{n}\delta_{u} + \bar{\sigma}_{\mathrm{odd},u}
\bigr)^{k},
\end{equation}
and the constant does not depend on $D_{\mathrm{odd}}$. Hence, as before,
in order to obtain an MA QC design, one needs to find $\bar{S}$ so as to
maximize $\bar{m}_{\mathrm{odd},3}$, then minimize
$\bar{m}_{\mathrm{odd},4}$, then maximize $\bar{m}_{\mathrm{odd},5}$, and
so on. In particular, analogously to (\ref{eq12}), $\bar{m}_{\mathrm{odd},3}=
\mathrm{constant} +F_{\mathrm{odd},3}$, where
\begin{equation}
\label{eq18} F_{\mathrm{odd},3} = 3 \bigl(2^{n} \bigr) \sum
_{u \in
\Delta_{0}}\bar{\sigma}_{\mathrm{odd},u}^{2}+ \sum
_{u \in\Delta} \bar{\sigma}_{\mathrm{odd},u}^{3}.
\end{equation}

We now have the following counterpart of Theorem \ref{th1} for an odd number of
factors.

\begin{theorem}\label{th3}
\textup{(a)} The inequality $F_{\mathrm{odd},3} \le
3(2^{2n})\{
2(v - s) - 1\}^{2}$ holds for every~$\bar{S}$.

\textup{(b)} Equality holds in \textup{(a)} if and only if $\bar{S}$
is even.
\end{theorem}

For $v - s \le2^{n - 1}$, by Theorem \ref{th3} and Lemma \ref{le2}, only even sets
$\bar{S} = \{ g_{s + 1},\ldots, g_{v}\}$ need to be considered in order to find
an MA QC design. Then the matrix $\bar{G} = [g_{s + 1}\cdots g_{v}]$ can be
represented as in (\ref{eq13}) via a binary matrix $B$ and the link
with the
associated regular design $d$ in $v - s - 1$ two-level factors and
$2^{n -
1}$ runs is again useful. As usual, denote the WLP of $d$ by
$A_{k}(d)$, $k
\ge1$, where $A{}_{1}(d) = A_{2}(d) = 0$. Also, write $A_{0}(d) = 1$ and
$A_{k}(d) = 0$ for $k > v - s - 1$. For $r = 2,3,\ldots,$ define
\begin{equation}
\label{eq19} E_{2r}(d)=\sum_{k = 0}^{2r}
\pmatrix{ %
v - s - 1 - k
\cr
\langle r - k/2 \rangle }
2^{k}A_{k}(d),
\end{equation}
where the combination is interpreted as zero if $\langle r - k/2 \rangle$, which is the
largest integer in $r - k/2$, exceeds $v - s - 1 - k$.

\begin{theorem}\label{th4}
Let $v - s \le2^{n - 1}$. Then
$\bar{S}$ yields an MA QC design in $2s + 1$ factors and
$2^{2n}$ runs if and only if the matrix B is so chosen that the
associated two-level regular design d sequentially minimizes $E_{2r}(d)$,
for $r = 2,3,\ldots,$ {etc}.
\end{theorem}

Since $\bar{\sigma}_{\mathrm{odd},u}$ and $\bar{m}_{\mathrm{odd},k}$ are
more involved than their counterparts in Section~\ref{sec2}, it is natural
that in
general the $E_{2r}(d)$ in Theorem \ref{th4} look more complicated than
the $A_{2r
- 1}(d) + A_{2r}(d)$ in Theorem \ref{th2}. By (\ref{eq19}), however,
$E_{4}(d)= \mathrm{constant}
+ 8A_{3}(d) + 16A_{4}(d)$, where the constant does not depend on $d$. So,
minimization of $E_{4}(d)$ simply calls for that of $A_{3}(d) + 2A_{4}(d)$
and, as shown below, this alone is often helpful.

\setcounter{example}{0}
\begin{example}[(continued)] Let $n = 4$, that is, $v = 120$, and $v - s= 5$.
Then out of the two nonisomorphic choices of $B$, namely, [1 2 12 3]
and [1
2 3 123], the first one uniquely minimizes $A_{3}(d) + 2A_{4}(d)$ and
hence, by Theorem \ref{th4}, yields an MA QC design in 231 factors and 256 runs.
\end{example}

\begin{example}[(continued)] Let $n = 5$, that is, $v = 496$, and $v - s= 10$.
Then among the five nonisomorphic choices of $B$ shown earlier, the one in
(ii) uniquely minimizes $A_{3}(d) + 2A_{4}(d)$ and hence, by Theorem \ref{th4},
yields an MA QC design in 973 factors and 1024 runs.
\end{example}

Indeed, for each $n = 3,4,5$ and $v - s \le2^{n - 1}$, one can check that
either (a)~there is a unique $B$ up to isomorphism or (b) the $B$ shown in
Table~\ref{tab1} uniquely minimizes $A_{3}(d) + 2A_{4}(d)$ and hence, by
Theorem \ref{th4},
leads to an MA QC design in $2s + 1$ factors and $2^{2n}$ runs. For $v
- s=
1$, the matrix $B$ does not arise in~(\ref{eq13}) and one only has to take
$\bar{S}$ as the singleton set consisting of $(1,0,\ldots,0)'$. Therefore, in
conjunction with what was found in Section~\ref{sec2}, we get, in particular, MA QC
designs in (i) 64 runs and 48--56 factors, (ii) 256 runs and 224--240
factors, and (iii)~1024 runs and 960--992 factors. While the MA designs in
(i) can be seen to agree with those reported in \cite{b16}, the ones in
(ii) and
(iii) are new.

\section{\texorpdfstring{Quaternary code designs in $2^{2n-1}$ runs}
{Quaternary code designs in 2 2n-1 runs}}\label{sec4}

Consider again the $2^{2n}$-run QC designs $D$
and $D_{\mathrm{odd}}$, in $2s$ and $2s + 1$ factors, introduced in
Sections~\ref{sec2} and \ref{sec3}. Suppose the last rows of the generator matrices, [$g_{1}\cdots g_{s}$]
for $D$ and [$g_{1}\cdots g_{s}\ g_{s + 1}$] for $D_{\mathrm{odd}}$, have all
elements even. Then by (\ref{eq5}), (\ref{eq6}) and (\ref{eq14}), the
row indexed by $u =
(u_{1},\ldots,u_{n})'$ in either design remains the same if $u_{n}$ is replaced
by $u_{n} + 2$ (mod 4), that is, the $2^{2n}$ runs of the design can be
split into two identical halves. Following \cite{b16}, any one of these halves
represents a $2^{2n-1}$-run QC design. It has the same number of factors
and, by (\ref{eq1}), the same WLP as the corresponding original design
$D$ or
$D_{\mathrm{odd}}$, and hence has MA if and only if so does the original
design. Therefore, depending on whether the number of factors is even or
odd, it suffices to find an MA QC design $D$ or $D_{\mathrm{odd}}$
based on
a generator matrix as stated above and then take one of its two identical
halves as the final design in $2^{2n-1}$ runs.

To adapt the complementary set theory for this purpose, observe that not
all vectors in the reference set $\Omega$ now qualify as columns of the
generator matrix, but that only the ones with last element even do so.
Write $\Omega_{0}=\{ g_{1},\ldots,g_{v_{0}}\}$ for the collection of these
qualifying vectors, where $v_{0} = 4^{n - 1} - 2^{n - 1}$. Let $\bar{S} =
\{ g_{s + 1},\ldots,g_{v_{0}}\}$ be the complement of $S = \{
g_{1},\ldots,g_{s}\}$ in $\Omega_{0}$. Define $\sigma_{u}$,
$\bar{\sigma}_{u}$, $m_{k}(D)$ and $\bar{m}_{k}$ as in Section~\ref{sec2}, and
$\sigma_{\mathrm{odd},u}$, $\bar{\sigma}_{\mathrm{odd},u}$,
$m_{k}(D_{\mathrm{odd}})$ and $\bar{m}_{\mathrm{odd},k}$ as in Section~\ref{sec3},
with $v$ replaced by $v_{0}$ in $\bar{\sigma}_{u}$ and
$\bar{\sigma}_{\mathrm{odd},u}$; cf. (\ref{eq9}) and (\ref{eq16}). Then
one can check that
$\sigma_{u}=2s$, $\bar{\sigma}_{u}=2(v_{0} - s)$,
$\sigma_{\mathrm{odd},u}=2s + 1$ and $\bar{\sigma}_{\mathrm{odd},u}=
2(v_{0} - s) - 1$ if $u$ has first $n - 1$ elements 0 and last element 0 or
2, and that the conclusion of Lemma \ref{le1} remains unaltered for every other
$u$. Therefore, despite working with $\Omega_{0}$ rather than $\Omega$, we
still have $m_{k}(D)= \mathrm{constant} +( - 1)^{k}\bar{m}_{k}$ and
$m_{k}(D_{\mathrm{odd}})= \mathrm{constant} + ( - 1)^{k}\bar{m}_{\mathrm{odd},k}$
for each $k$, where the constants do not depend on $D$ or
$D_{\mathrm{odd}}$. Furthermore, in the representation (\ref{eq13}) for
an even set
$\bar{S}$ via $\bar{G}$, each vector in $\bar{S}$ has last element even and
hence $\bar{S} \subset\Omega_{0}$, as required here.

From the above, it is evident that the findings in Sections~\ref{sec2} and \ref{sec3} as well
as Table~\ref{tab1} continue to remain valid in the present setup, with $v$ simply
replaced by $v_{0}$, thus leading to $2^{2n-1}$-run MA QC designs in $2s$
factors, $0 \le v_{0} - s \le2^{n - 1}$, and $2s + 1$ factors, $1 \le v_{0}
- s \le2^{n - 1}$. As a result, we get, in particular, MA QC designs in
(i) 128 runs and 96--112 factors, and (ii) 512 runs and 448--480 factors.

\begin{example}\label{ex3}
Let $n = 4$, that is, $v_{0} = 56$, and $v_{0} - s = 6$. Then
from Table~\ref{tab1}, $B = {}$[1 2 12 3 13] and, hence, by (\ref{eq13}),
\[
\bar{S}= \bigl\{ (1,0,0,0)',(1,2,0,0)',(1,0,2,0)',(1,2,2,0)',(1,0,0,2)',(1,2,0,2)'
\bigr\}.
\]
If one (i) finds the complement $S = \{ g_{1},\ldots,g_{s}\}$ of $\bar{S}$ in
$\Omega_{0}$ and (ii) constructs a design $D$ in 100 factors and 256 runs
from the generator matrix [$g_{1}\cdots g_{s}$], then any one of the two
identical halves of $D$ is an MA QC design in 100 factors and 128 runs.
\end{example}

\begin{example}\label{ex4}
Let $n = 5$, that is, $v_{0} = 240$, and $v_{0} - s= 7$. Then
from Table~\ref{tab1}, $B = {}$[1 2 12 3 4 34] and, hence, by (\ref{eq13}),
\begin{eqnarray*}
\bar{S}&=& \bigl\{ (1,0,0,0,0)',(1,2,0,0,0)',(1,0,2,0,0)',(1,2,2,0,0)',
(1,0,0,2,0)',
\\
&&\hspace*{200pt}(1,0,0,0,2)', (1,0,0,2,2)' \bigr\}.
\end{eqnarray*}
Now, if one (i) finds the complement $S = \{ g_{1},\ldots,g_{s}\}$ of
$\bar{S}$ in $\Omega_{0}$, (ii) constructs a design $D_{\mathrm{even}}$ in
468 factors and 1024 runs from the generator matrix [$g_{1}\cdots g_{s}\ g_{s +
1}$], where $g_{s + 1}= (1,0,0,0,0)'$, and (iii) deletes the last
column of
$D_{\mathrm{even}}$ (i.e., the second column contributed by $g_{s +
1}$) to
get a design $D_{\mathrm{odd}}$ in 467 factors and 1024 runs, then any one
of the two identical halves of $D_{\mathrm{odd}}$ is an MA QC design in 467
factors and 512 runs.
\end{example}

\section{Comparison with regular MA designs and concluding
remarks}\label{sec5} In a highly fractionated setup, our results
explore the best that can be achieved by QC designs under the MA
criterion and hence facilitate comparison with their regular
counterparts. Indeed, as seen below, MA QC designs obtained here
compete very well with MA regular designs. For illustration, we
consider the cases of $N = 128$ and 256 runs and recall that for these
$N$, our results yield MA QC designs for $96 \le q \le112$ and $224 \le
q \le240$, where $q$ is the number of factors.

For $N = 128$, MA QC designs have (i) the same WLP as MA regular
designs if\vadjust{\goodbreak}
$96 \le q \le99$ or $109 \le q \le112$, and (ii) the same $A_{3}$ but a
little larger $A_{4}$ if $100 \le q \le108$. In the first case, MA QC
designs clearly outperform MA regular designs because of higher resolution
and projectivity. For most practical purposes, in the second case too the
same features of MA QC designs far outweigh their slightly higher
aberration. For instance, if $q = 103$, then the $A_{4}$ values for the MA
QC and MA regular designs are 35,707 and 35,705, respectively, while both
have $A_{3}= 1360$. Thus, the marginally larger $A_{4}$ for the MA QC
design is more than compensated by the fact that it has projectivity at
least three, while the MA regular design has as many as 1360 three-factor
projections which do not contain a complete $2^{3}$ factorial.

For $N = 256$, MA QC designs have (i) the same WLP as MA regular
designs if
$224 \le q \le227$ or $237 \le q \le240$, and (ii) the same $A_{3}$ but
marginally larger $A_{4}$ if $228 \le q \le236$, for example, if $q =
228$, then the $A_{4}$ values for the MA QC and MA regular designs are 434,057
and 434,056, respectively, while both have $A_{3}= 7616$. This has the same
implications as before in favor of MA QC designs. The same pattern is seen
to persist for larger $N$.

It will be of interest to extend the present results to QC designs which
are less highly fractionated than the ones considered here, that is, for
which the size of $\bar{S}$ exceeds $2^{n - 1}$. In view of Lemma \ref{le2}, then
the bounds in Theorems \ref{th1} and \ref{th3} are not attainable and, therefore, one
cannot have neat results in terms of even sets. However, Lemma \ref{le1} as
well as
Lemma \ref{lea1} and equations (\ref{eqA.4}), (\ref{eqA.5}), (\ref{eqA.15}),
(\ref{eqA.16}) in the
\hyperref[app]{Appendix}
continue to hold and should be useful in replacing the bounds in
Theorems \ref{th1}
and \ref{th3} by sharper, attainable ones. We conclude with the hope that the
present endeavor will generate more interest in this area.

\begin{appendix}\label{app}
\section*{Appendix: Proofs}
\begin{pf*}{Proof of Lemma \ref{le1}}
For $u \in\Delta$, consider the row of [$D\bar{D}$]
which is indexed by $u$. This row is of the form (\ref{eq6}) with $s$ there
replaced by $v$. Hence, by (\ref{eq5}), analogously to (\ref{eq9}),
\begin{equation}
\label{eqA.1} %(A$.1)
\sigma_{u} + \bar{\sigma}_{u}=
\sum_{j = 1}^{v} \bigl(i^{u'g_{j}} +
i{}^{ - u'g_{j}} \bigr).
\end{equation}

By the definitions of $\Omega$, $\Delta$ and $\Delta_{0}$, the union of the
sets $\{ g_{j}, - g_{j}\}$, $1 \le j \le v$, equals $\Delta- \Delta_{0}$.
Hence by (\ref{eqA.1}), $\sigma_{u} + \bar{\sigma}_{u}=\sum_{w \in
\Delta}
i^{u'w} - \sum_{w \in\Delta_{0}}i^{u'w}$. The lemma now follows, noting
that (i) $\sum_{w \in\Delta} i^{u'w}= 0$ for nonnull $u$, (ii)
$\sum_{w \in\Delta_{0}}i^{u'w}= 0$ if $u$ has an odd element, and (iii)
$\sum_{w \in\Delta_{0}}i^{u'w} = 2^{n}$ if $u$ has all elements even,
since then $u'w= 0$ mod 4, for every $w \in\Delta_{0}$.
\end{pf*}

In order to prove Theorem \ref{th1}, we require some notation and another lemma.
For any $n \times1$ vector $g$ with integer elements, define $\alpha(g)$
as 1 or 0 according to whether $g$ is null (mod 4) or not. Recall
that\vadjust{\goodbreak}
$\bar{S} = \{ g_{s + 1},\ldots,g_{v}\}$. For $s + 1 \le j,k,h \le v$, let
\begin{eqnarray}
\label{eqA.2} %(\mathrm{A}.2)
\beta_{jkh}&=&\alpha(g_{j} +
g_{k} + g_{h})+\alpha(g_{j} + g_{k} -
g_{h})+ \alpha(g_{j} - g_{k} + g_{h})
\nonumber
\\[-8pt]
\\[-8pt]
\nonumber
&&{}+ \alpha(g_{j} - g_{k} - g_{h}).
\end{eqnarray}

\begin{lemma}\label{lea1}
\textup{(a)} For each $j,k,h$, $\beta_{jkh} = 0$ or 1.

\textup{(b)} For any fixed $j,k$, \textup{(i)} $\beta_{jkh} = 0$ for all $h$,
if $g_{j} + g_{k}$ has all elements even, \textup{(ii)}~$\beta_{jkh} = 1$
for at most two choices of h, otherwise.
\end{lemma}

\begin{pf}
(a) This follows noting that the right-hand side of (\ref
{eqA.2}) cannot
have a pair of terms both of which equal 1. For instance, if the first two
terms equal 1, then $g_{j} + g_{k} + g_{h} = 0$ (mod 4) and $g_{j} + g_{k}
- g_{h} = 0$ (mod 4). So $2g_{h}= 0$ (mod~4), which is impossible as $g_{h}$
has an odd element. The same argument applies to any other pair of terms.

(b) If $g_{j} + g_{k}$ has all elements even, then the same holds for
$g_{j} - g_{k}$. As $g_{h}$ has an odd element, then all terms on the
right-hand side of (\ref{eqA.2}) vanish, that is, $\beta_{jkh}= 0$.
Next, let
$g_{j} + g_{k}$ have an odd element. Then, by (\ref{eqA.2}), $\beta
_{jkh}= 1$ if
and only if $g_{h} = \pm(g_{j} + g_{k})$ (mod 4) or $g_{h}= \pm(g_{j} -
g_{k})$ (mod 4). Since $- (g_{j} + g_{k}) = 3(g_{j} + g_{k})$ (mod 4) and $-
(g_{j} - g_{k}) = 3(g_{j} - g_{k})$ (mod 4), and no two vectors in
$\bar{S}$ can be multiples of each other, it follows that $\beta_{jkh}= 1$
for at most two choices of $h$.
\end{pf}

\begin{pf*}{Proof of Theorem \ref{th1}} (a) As in the proof of Lemma \ref{le1} but using a more formal
notation,
\begin{equation}
\label{eqA.3} %(\mathrm{A}.3)
\qquad\sum_{u \in\Delta_{0}}i^{u'g}
= \sum_{u \in
\Delta_{0}}i^{ - u'g}= 2^{n}
\alpha(2g),\qquad \sum_{u \in\Delta} i^{u'g} = \sum
_{u \in\Delta} i^{ - u'g}=2^{2n}\alpha(g).
\end{equation}
Let $\sum^{(2)}$ and $\sum^{(3)}$ denote double and triple sums on $j,k$
and $j,k,h$ over the ranges $s + 1 \le j,k \le v$ and $s + 1 \le j,k,h
\le
v$, respectively. By (\ref{eq9}) and (\ref{eqA.3}),
\begin{eqnarray}
\label{eqA.4} %(\mathrm{A}.4)
\sum_{u \in\Delta_{0}}\bar{
\sigma}_{u}^{2}&=& {\sum}^{(2)}\sum
_{u \in\Delta_{0}} \bigl(i^{u'g_{j}} + i^{ -
u'g_{j}}
\bigr) \bigl(i^{u'g_{k}} + i^{ - u'g_{k}} \bigr)
\nonumber
\\[-8pt]
\\[-8pt]
\nonumber
&=&2^{n + 2}{\sum}^{(2)}
\alpha(2g_{j} + 2g_{k}),
\end{eqnarray}
because $2g_{j} - 2g_{k}=2g_{j} + 2g_{k}$ (mod 4). Similarly,
\begin{eqnarray}
\label{eqA.5} %(A$.5)
\sum_{u \in\Delta} \bar{
\sigma}_{u}^{3}&=& {\sum}^{(3)}
\sum_{u
\in
\Delta} \bigl(i^{u'g_{j}} + i^{ - u'g_{j}}
\bigr) \bigl(i^{u'g_{k}} + i^{ -
u'g_{k}} \bigr) \bigl(i^{u'g_{h}} +
i^{ - u'g_{h}} \bigr)
\nonumber
\\[-8pt]
\\[-8pt]
\nonumber
&=&2^{2n + 1} {\sum}^{(3)}\beta
{}_{jkh},
\end{eqnarray}
where $\beta_{jkh}$ is given by (\ref{eqA.2}).

Observe that $\bar{S} = \{ g_{s + 1},\ldots,g_{v}\}$ can be partitioned into
$t\ ( \ge1)$ mutually exclusive and exhaustive nonempty subsets such that
any $g_{j} + g_{k}$, $s + 1 \le j,k \le v$, has all elements even if and
only if $g_{j}$ and $g_{k}$ belong to the same subset. Denote the
cardinalities of these subsets by $f_{1},\ldots,f_{t}$. Then
\begin{equation}
\label{eqA.6} %(A$.6)
f_{1} +\cdots + f_{t} = v - s,
\end{equation}
since $\bar{S}$ has cardinality $v - s$. Clearly, $\bar{S}$ is even if and
only if $t = 1$. Now by~(\ref{eqA.4}),
\begin{equation}
\label{eqA.7} %(A$.7)
\sum_{u \in\Delta_{0}}\bar{
\sigma}_{u}^{2}=2^{n + 2} \bigl(f_{1}^{2}
+\cdots + f_{t}^{2} \bigr),
\end{equation}
because $\alpha(2g_{j} + 2g_{k}) = 1$ if and only if $g_{j} + g_{k}$ has
all elements even, that is, $g_{j}$ and $g_{k}$ belong to the same
subset of
$\bar{S}$ as described above. Similarly, by (\ref{eqA.5}), (\ref
{eqA.6}) and Lemma \ref{lea1},
\begin{equation}
\label{eqA.8} %(\mathrm{A}.8)
\sum_{u \in\Delta} \bar{
\sigma}_{u}^{3} \le2^{2n +
2}\sum
_{l = 1}^{t}\sum_{r( \ne l) = 1}^{t}f_{l}f_{r}=2^{2n + 2}
\bigl\{ (v - s)^{2} - \bigl(f_{1}^{2} +\cdots +
f_{t}^{2} \bigr) \bigr\}.
\end{equation}
By (\ref{eq12}) and (\ref{eqA.6})--(\ref{eqA.8}),
\begin{equation}
\label{eqA.9} %(A$.9)
F_{3} \le2^{2n + 2} \bigl\{ (v -
s)^{2} + 2 \bigl(f_{1}^{2} +\cdots +
f_{t}^{2} \bigr) \bigr\} \le3 \bigl(2^{2n + 2} \bigr)
(v - s)^{2},
\end{equation}
which proves (a).

(b) By (\ref{eqA.6}) and (\ref{eqA.9}), equality holds in (a) only if
$f_{1}^{2} +\cdots +
f_{t}^{2}=(f_{1} +\cdots + f_{t})^{2}$, which holds only if $t =1$, that is,
the set $\bar{S}$ is even. On the other hand, if $\bar{S}$ is even,
then by
(\ref{eqA.4}), (\ref{eqA.5}) and Lemma \ref{lea1},
\[
\sum_{u \in\Delta_{0}}\bar{\sigma}_{u}^{2}=2^{n + 2}(v
- s)^{2}, \qquad \sum_{u \in\Delta} \bar{
\sigma}_{u}^{3}= 0.
\]
Therefore, (\ref{eq12}) yields $F_{3} = 3(2^{2n + 2})(v - s)^{2}$, and equality
holds in (a).
\end{pf*}

\begin{pf*}{Proof of Lemma \ref{le2}} Only if: For an even set $\bar{S} = \{ g_{s +
1},\ldots,g_{v}\}$, the vectors $g_{s + 1} + g_{j}$ (mod 4), $s + 1 \le j
\le
v$, satisfy the following: (i) each of them has all elements even, and (ii)
no two of them add up to $2g_{s + 1}$ (mod~4). Here (ii) is due to the fact
that no two vectors in $\bar{S}$ are multiples of each other over $Z_{4}$.
Since there are at most $2^{n - 1}$ distinct $n \times1$ vectors over
$Z_{4}$ satisfying (i) and (ii), the only if part follows.

If: There are $2^{n - 1}$ distinct $(n - 1) \times1$ vectors over $Z_{4}$,
each of which has all elements even. For $v - s \le2^{n - 1}$, consider
any $v - s$ of these $(n - 1) \times1$ vectors, say, $\tilde{g}_{s +
1},\ldots,\tilde{g}_{v}$. Then the set $\bar{S}$, consisting of the vectors
$(1,\tilde{g}'_{j})'$, $s+1 \le j \le v$, is even and the if part
follows.
\end{pf*}

\begin{pf*}{Proof Of Theorem \ref{th2}} Let $b_{s + 2},\ldots,b_{v}$ denote the columns of $B$.
For any $u = (u_{1},\ldots,u_{n})' \in\Delta$, writing $u(2) =
(u_{2},\ldots,u_{n})'$, from (\ref{eq9}) and (\ref{eq13}),
\begin{equation}
\label{eqA.10} \bar{\sigma}_{u}= \bigl(i^{u_{1}} +
i^{ - u_{1}} \bigr) \Biggl\{ 1 + \sum_{j = s +
2}^{v}(
- 1)^{u(2)'b_{j}} \Biggr\}.
\end{equation}
Therefore, considering $u_{1} = 0,1,2,3$ separately, for any positive
integer $k$,
\begin{equation}
\label{eqA.11} %(A$.11)
\sum_{u \in\Delta} \bar{
\sigma}_{u}^{k}=2^{k} \bigl\{ 1 + ( -
1)^{k} \bigr\} \sum_{u_{2} = 0}^{3}
\cdots \sum_{u_{n} = 0}^{3} \Biggl\{ 1 + \sum
_{j = s +
2}^{v}( - 1)^{u(2)'b_{j}} \Biggr
\}^{k}.
\end{equation}
Note that the quantities $( - 1)^{u(2)'b_{j}}$, $s + 2 \le j \le v$, remain
unaltered if any element $u_{h}$ of $u(2)$ is replaced by $u_{h} + 2$ (mod
4). Hence, (\ref{eqA.11}) yields
\begin{eqnarray}
\label{eqA.12} %(\mathrm{A}.12)
\sum_{u \in\Delta} \bar{
\sigma}_{u}^{k}&=&2^{n + k -
1} \bigl\{ 1 + ( -
1)^{k} \bigr\} \sum_{u_{2} = 0}^{1}
\cdots \sum_{u_{n} = 0}^{1} \Biggl\{ 1 + \sum
_{j = s + 2}^{v}( - 1)^{u(2)'b_{j}} \Biggr
\}^{k}
\nonumber
\\[-8pt]
\\[-8pt]
\nonumber
&=&2^{n + k - 1} \bigl\{ 1 + ( - 1)^{k} \bigr\} \sum
_{x \in\Gamma} (1 + \lambda_{x})^{k},
\end{eqnarray}
where $\Gamma$ is the collection of the $2^{n - 1}$ binary column vectors
of order $n - 1$, and for any $x \in\Gamma$,
\begin{equation}
\label{eqA.13} %(A$.13)
\lambda_{x} = \sum
_{j = s + 2}^{v}( - 1)^{x'b_{j}}.
\end{equation}
Also, by (\ref{eqA.10}), for every $u = (u_{1},\ldots,u_{n})' \in\Delta_{0}$,
$\bar{\sigma}_{u}$ equals $2(v - s)$ or $- 2(v - s)$, according to whether
$u_{1}= 0$ or 2, respectively. For any nonnegative integer $r$, therefore,
$\sum_{u \in\Delta} \delta_{u}\bar{\sigma}_{u}^{r}$ is a constant which
does not depend on the choice of $\bar{S}$. Hence, by (\ref{eq11}) and
(\ref{eqA.12}),
\begin{eqnarray*}
\bar{m}_{k}&=& \mathrm{constant} + 2^{n + k} \sum
_{x \in\Gamma} (1 + \lambda_{x})^{k}\qquad
\mbox{if $k$ is even,}
\\
&=& \mathrm{constant}\qquad \mbox{if $k$ is odd,}
\end{eqnarray*}
where the constants do not depend on $\bar{S}$. Therefore, $\bar{S}$ yields
an MA QC design if and only if the matrix $B$ is so chosen as to
sequentially minimize $\sum_{x \in\Gamma} (1 + \lambda_{x})^{2r}$, for
$r = 2,3,\ldots,$ etc.

Now from (\ref{eqA.13}), observe that $\lambda_{x}$, $x \in\Gamma$,
are the row
(run) sums of the two-level regular design $d$ generated by $B = [b_{s +
2}\cdots b_{v}]$. Recall that $d$ involves $v - s - 1$ two-level factors and
$2^{n - 1}$ runs. Write $\tilde{d}$ for the full foldover of $d$ and note
the following:\vadjust{\goodbreak}
\begin{longlist}[(iii)]
\item[(i)] $\tilde{d}$ involves $v - s$ two-level factors and $2^{n}$ runs,

\item[(ii)] the row (run) sums of $\tilde{d}$ are $\pm(1 + \lambda_{x})$, $x
\in
\Gamma$, because those of $d$ are $\lambda_{x}$, $x \in\Gamma$,

\item[(iii)] in any two-level regular $N$-run design with rows (runs) $\xi
'_{1},\ldots,\xi'_{N}$, the $N$ row sums occur equally often among the
$N^{2}$ scalar products $\{\xi'_{j}\xi_{k}\dvtx1 \le j,k \le N\}$ of the
rows.
\end{longlist}

So, as in the passage from (\ref{eq4}) to (\ref{eq8}) in Section~\ref{sec2},
sequential minimization
of $\sum_{x \in\Gamma} (1 + \lambda_{x})^{2r}$, for $r = 2,3,\ldots,$ etc.
amounts to choosing $d$ so that its foldover $\tilde{d}$ has MMA, and hence
MA, among all such foldovers. Now the result follows if we denote the WLP
of $\tilde{d}$ by $A_{k}(\tilde{d})$, $k \ge1$, and note that
$A_{1}(\tilde{d}) = A_{2}(\tilde{d}) = 0$, while $A_{2r - 1}(\tilde{d}) =
0$, $A_{2r}(\tilde{d}) = A_{2r - 1}(d) + A_{2r}(d)$, for $r = 2,3,\ldots,$ etc.
\end{pf*}

We indicate only the key steps in the proofs of Theorems \ref{th3} and \ref{th4} which are
similar to but more elaborate than those of Theorems \ref{th1} and \ref{th2}.

\begin{pf*}{Proof of Theorem \ref{th3}}
(a) By (\ref{eq5}) and (\ref{eq16}), analogously to
(\ref{eq9}),
\begin{equation}
\label{eqA.14} %(\mathrm{A}.14)
\qquad\bar{\sigma}_{\mathrm{odd},u}=
\bigl(i^{u'g_{s + 1}} + i^{3 -
u'g_{s + 1}} \bigr)/(1 - i)+\sum
_{j = s + 2}^{v} \bigl(i^{u'g_{j}} + i^{ - u'g_{j}}
\bigr), \qquad u \in\Delta.
\end{equation}
Thus, using (\ref{eqA.3}), along the lines of (\ref{eqA.4}) and (\ref
{eqA.5}), but with heavier
algebra,
\begin{eqnarray}
\label{eqA.15}\qquad  %(A$.15)
\sum_{u \in\Delta_{0}}\bar{
\sigma}_{\mathrm
{odd},u}^{2}&=&2^{n} + 2^{n + 2} \Bigl\{
\tilde{\sum}^{(1)}\alpha(2g_{s + 1} + 2g_{j})
+ \tilde{\sum}^{(2)}\alpha(2g_{j} +
2g_{k}) \Bigr\},
\\
%(A$.16)
\label{eqA.16} \sum_{u \in\Delta} \bar{
\sigma}_{\mathrm{odd},u}^{3}&=& 3 \bigl(2^{2n} \bigr)\tilde{\sum
}^{(2)}\beta{}_{s
+ 1jk}+2^{2n + 1} \tilde{\sum
}^{(3)}\beta{}_{jkh},
\end{eqnarray}
where $\tilde{\sum}^{(1)}$ denotes the sum on $j$ over $s + 2 \le j
\le
v$, while $\tilde{\sum}^{(2)}$ and $\tilde{\sum}^{(3)}$ denote double
and triple sums on $j,k$ and $j,k,h$ over $s + 2 \le j,k \le v$ and $s
+ 2
\le j,k,h \le v$, respectively.

Partition $\bar{S} = \{ g_{s + 1},\ldots,g_{v}\}$ into $t\ ( \ge1)$ mutually
exclusive and exhaustive nonempty subsets as specified in the proof of
Theorem \ref{th1}. The cardinalities $f_{1},\ldots,f_{t}$ of these subsets satisfy
(\ref{eqA.6}). Without loss of generality, let $g_{s + 1}$ belong to
the first of
these subsets. Then $\tilde{\sum}^{(1)}\alpha(2g_{s + 1} +
2g_{j})=f_{1}^{*}$, while, analogously to~(\ref{eqA.7}),
$\tilde{\sum}^{(2)}\alpha(2g_{j} + 2g_{k})= \sum_{l =
1}^{t}(f_{l}^{*})^{2}$, where $f_{1}^{*} = f_{1} - 1$ and $f_{l}^{*} =
f_{l}$ for $l \ge2$. Hence, by (\ref{eqA.15}),
\begin{equation}
\label{eqA.17} %(A$.17)
\sum_{u \in\Delta_{0}}\bar{
\sigma}_{\mathrm
{odd},u}^{2}=2^{n} \Biggl\{ 1 +
4f_{1}^{*} + 4\sum_{l = 1}^{t}
\bigl(f_{l}^{*} \bigr)^{2} \Biggr\}.
\end{equation}
Also, by Lemma \ref{lea1}, $\tilde{\sum}^{(2)}\beta{}_{s + 1jk} \le2\sum_{l
= 2}^{t}f_{l}^{*}$, while, analogously to (\ref{eqA.8}), $\tilde{\sum }^{(3)}\beta
{}_{jkh} \le2\{ (\sum_{l = 1}^{t}f_{l}^{*})^{2} - \sum_{l =
1}^{t}(f_{l}^{*})^{2}\}$, so that by (\ref{eqA.16}),
\begin{equation}
\label{eqA.18} %(A$.18)
\sum_{u \in\Delta} \bar{
\sigma}_{\mathrm{odd},u}^{3} \le 2^{2n +
1} \Biggl\{ 3\sum
_{l = 2}^{t}f_{l}^{*} + 2
\Biggl(\sum_{l = 1}^{t}f_{l}^{*}
\Biggr)^{2} - 2\sum_{l = 1}^{t}
\bigl(f_{l}^{*} \bigr)^{2} \Biggr\}.
\end{equation}
By (\ref{eqA.6}), the sum of the nonnegative integers $f_{1}^{*},\ldots,f_{t}^{*}$
equals $v - s - 1$. Therefore,
\[
f_{1}^{*} \le v - s - 1,\qquad \sum
_{l = 1}^{t} \bigl(f_{l}^{*}
\bigr)^{2} \le(v - s - 1)^{2},
\]
and, hence, from (\ref{eq18}), (\ref{eqA.17}) and (\ref{eqA.18}), on
simplification
\begin{eqnarray*}
F_{\mathrm{odd},3} &\le&2^{2n} \Biggl\{ 3 + 6(v - s - 1) + 4(v - s -
1)^{2} + 6f_{1}^{*} + 8\sum
_{l = 1}^{t} \bigl(f_{l}^{*}
\bigr)^{2} \Biggr\} \\
&\le&3 \bigl(2^{2n} \bigr) \bigl\{ 2(v - s) -
1 \bigr\}^{2},
\end{eqnarray*}
which proves (a).

(b) It is easily seen that equality holds in (a) only if $f_{1}^{*} = v
- s
- 1$, that is, $f_{1} = v - s$, which holds only if $\bar{S}$ is even. On
the other hand, if $\bar{S}$ is even, then by~(\ref{eqA.15}), (\ref
{eqA.16}) and Lemma~\ref{lea1},
\[
\sum_{u \in\Delta_{0}}\bar{\sigma}_{\mathrm{odd},u}^{2}=2^{n}
+ 2^{n +
2} \bigl\{ v - s - 1 + (v - s - 1)^{2} \bigr
\}=2^{n} \bigl\{ 2(v - s) - 1 \bigr\}^{2},
\]
and $\sum_{u \in\Delta} \bar{\sigma}_{\mathrm{odd},u}^{3}= 0$, so that
by (\ref{eq18}), equality holds in (a).
\end{pf*}

Some notation and a lemma are needed for proving Theorem \ref{th4}. With $B$ as in~(\ref{eq13}), the binary matrix [$B\ B$] generates a regular design, say,
$d_{0}$, in
$2(v - s - 1)$ two-level factors and $2^{n - 1}$ runs. Denote the full
foldover of $d_{0}$ by $\tilde{d}_{0}$. Let $A_{k}(d_{0})$ and
$A_{k}(\tilde{d}_{0})$, $k \ge1$, be the WLPs of $d_{0}$ and
$\tilde{d}_{0}$, respectively. Clearly, $A_{1}(d_{0}) = 0$ and
$A_{2}(d_{0}) = v - s - 1$, since the columns of [$B\ B$] are nonnull but
identical in pairs.

\begin{lemma}\label{lea2}
\textup{(a)} $A_{k}(\tilde{d}_{0}) = 0$, for every odd $k$, \textup{(b)}
$A_{2}(\tilde{d}_{0}) = v - s - 1$, \textup{(c)}~$A_{2r}(\tilde{d}_{0}) =
E_{2r}(d)$, for $r = 2,3,\ldots,$ where $E_{2r}(d)$
is given by (\ref{eq19}).
\end{lemma}

\begin{pf}
While (a) holds for any full foldover design, (b) is obvious. To
prove~(c), denote the columns of $B$ by $b_{s + 2},\ldots,b_{v}$ and write
$[B\ B]=[b_{s + 2}^{(1)}\cdots  b_{v}^{(1)}\ b_{s + 2}^{(2)}\cdots  b_{v}^{(2)}]$, where
$b_{j}^{(1)} = b_{j}^{(2)} = b_{j}$, $s + 2 \le j \le v$. Then any set of
$h$ columns of [$B\ B$], forming a word of length $h$ of the design $d_{0}$,
has the structure $\{
b_{j_{1}}^{(l_{1})},\ldots,b_{j_{k}}^{(l_{k})}\}\cup\{
b_{j}^{(1)},b_{j}^{(2)}\dvtx j \in J\}$, where $k$ is such that $h - k$
is a
nonnegative even integer, the columns $b_{j_{1}},\ldots,b_{j_{k}}$ constitute
a word of length $k$ of $d$, each of $l_{1},\ldots,l_{k}$ is either 1 or 2,
and $J$ is any subset, with cardinality $(h - k)/2$, of the complement of
$\{ j_{1},\ldots,j_{k}\}$ in $\{ s + 2,\ldots,v\}$. So
\begin{equation}
\label{eqA.19} %(A$.19)
A_{h}(d_{0})= \sum
\pmatrix{ v - s - 1 - k
\cr
(h - k)/2 }2^{k}A_{k}(d),
\end{equation}
where the sum ranges over $k = h,h - 2,\ldots,$ etc. Now,
$A_{2r}(\tilde{d}_{0}) = A_{2r - 1}(d_{0}) + A_{2r}(d_{0})$, since
$\tilde{d}_{0}$ is the full foldover of $d_{0}$. Hence, (c) follows from
(\ref{eqA.19}), recalling the definition of $E_{2r}(d)$ from (\ref
{eq19}).
\end{pf}

\begin{pf*}{Proof of Theorem \ref{th4}} By (\ref{eq13}) and
(\ref{eqA.14}), using the same
notation as in
(\ref{eqA.10}),
\[
\bar{\sigma}_{\mathrm{odd},u}= \bigl(i^{u_{1}} + i^{3 - u_{1}} \bigr)/(1
- i)+ \bigl(i^{u_{1}} + i^{ - u_{1}} \bigr)\sum
_{j = s + 2}^{v}( - 1)^{u(2)'b_{j}}, \qquad u \in\Delta.
\]

Hence, by (\ref{eq17}), arguing as in the proof of Theorem \ref{th2}, for any positive
integer $k$,
\begin{eqnarray*}
\bar{m}_{\mathrm{odd},k}&=& \mathrm{constant} + 2^{n}\sum
_{x \in\Gamma} (1 + 2\lambda_{x})^{k}\qquad
\mbox{if $k$ is even,}
\\
&=& \mathrm{constant}\qquad \mbox{if $k$ is odd,}
\end{eqnarray*}
where $\lambda_{x}$ is given by (\ref{eqA.13}) and the constants do not
depend on
$\bar{S}$. Therefore, $\bar{S}$ yields an MA QC design if and only if the
matrix $B$ is so chosen as to sequentially minimize $\sum_{x \in
\Gamma}
(1 + 2\lambda_{x})^{2r}$, for $r = 2,3,\ldots,$ etc. Again, as with Theorem
\ref{th2}, this happens if and only if the full foldover design $\tilde{d}_{0}$ in
Lemma \ref{lea2} has MA among all such foldovers, because the row (run) sums of
$\tilde{d}_{0}$ are $\pm(1 + 2\lambda_{x})$, $x \in\Gamma$, as those of
$d$ are $\lambda_{x}$, $x \in\Gamma$. The result is now immediate from
Lemma \ref{lea2}.
\end{pf*}
\end{appendix}
% imsref loaded by akundreckaite, 2013-11-11 15:43:51
%

% zodis "Acknowledgments" paliekamas pagal autoriu

%suskaldyti doi

\printaddresses

\end{document}